\documentclass[12pt]{amsart}
\usepackage{amssymb}

\newcommand{\bea}{\begin{eqnarray}}
\newcommand{\eea}{\end{eqnarray}}

\newtheorem{Thm}{Theorem}[section]

\newtheorem{Cor}[Thm]{Corollary}
\newtheorem{Lem}[Thm]{Lemma}
\newtheorem{Prop}[Thm]{Proposition}

\newcommand{\bd}{\begin{displaystyle}}
\newcommand{\ed}{\end{displaystyle}}

\newcommand{\sun}{\sum_{n \in Z}}

\newcommand{\suz}{{ \sum_{k \in \mathbb  Z}}}
\newcommand{\ltr}{L^{2}(\mathbb R)}

\def \R{\rm {\bf R}}

\def \Bbb{\mathbb}

\begin{document}

\title{Perturbations of Weyl-Heisenberg frames}
\author{Peter G. Casazza, Ole Christensen and Mark C. Lammers}
\address{Department of Mathematics \\
University of Missouri-Columbia \\
Columbia, MO 65211 and
Department of Mathematics, \\
Technical University of Denmark \\
 2800 Lyngby, Denmark, and
Department of Mathematics, University \\
of South Carolina, Columbia, SC 29208}

\thanks{The first author was supported by NSF DMS 970618}

\email{pete@math.missouri.edu;
Ole.Christensen@mat.dtu.dk; lammers @math.sc.edu}


\maketitle

\begin{abstract}
We develop a usable perturbation theory for Weyl-Heisenberg
frames.  In particular, we prove that if $(E_{mb}T_{na}g)_{m,n\in
\mathbb Z}$
is a WH-frame and $h$ is a function which is close to $g$ in the Wiener 
Amalgam
space norm, then also  $(E_{mb}T_{na}h)_{m,n\in \mathbb Z}$ is a
WH-frame.
We also prove perturbation results for the parameters $a,b$. 
\end{abstract}

\section{Introduction} \label{intro}

In 1952, Duffin and Schaeffer \cite{DS} introduced the notion of a frame
for a Hilbert space.  A sequence $(f_{n})_{n\in I}$ is a {\bf frame}
for a Hilbert space $H$ if there are constants $A,B>0$ satisfying,
\begin{equation}
A\|f\|^{2}\le \sum_{n\in I}|\langle f,f_{n} \rangle |^{2} \le
B\|f\|^{2},\label{frameident}
\end{equation}
for all $f\in H$.  The constant $A$ (respectively, $B$) is a {\bf lower}
(resp. {\bf upper}) frame bound for the frame.  One of the most
important 
frames
for applications, especially signal processing,
are the Weyl-Heisenberg frames.  For $g\in L^{2}(\mathbb R )$ we define
the
{\bf translation parameter} $a>0$ and the {\bf modulation parameter}
$b>0$
by:
$$
E_{mb}g(t) = e^{2{\pi}imbt},\ \ \ T_{na}g(t) = g(t-na).
$$
For $g\in L^{2}(\mathbb R )$ and $a,b>0$,
we say for short that $(g,a,b)$ {\bf is a Weyl-Heisenberg frame for}
$L^{2}(\mathbb R )$ if $(E_{mb}T_{na}g)_{m,n\in \mathbb Z}$ is a frame
for
$L^{2}(\mathbb R )$. 
We call $(f_{n})_{n\in I}$ a {\bf Riesz basis} (resp. {\bf Riesz
basic sequence}) for a Hilbert space $H$ if it is a bounded
unconditional basis for $H$ (resp. for its closed linear span.) \\ 

 Weyl-Heisenberg frames are
extremely sensitive to even arbitrarily small changes in the function
$g$ and the translation and modulation parameters. For example,
$(E_{m}T_{n}\chi_{[0,1]})_{m,n\in \mathbb Z}$ is a frame
for $L^{2}(\mathbb R )$, but for arbitrary $\epsilon >0$, the
functions $(E_{m}T_{n}\chi_{[0,1- \epsilon]})_{m,n\in \mathbb Z}$
are not. As a result, there
are
few general theorems on perturbations of Weyl-Heisenberg frames and
those
that exist are often very technical in nature (see \cite{FZ}, and the
article of Christensen in \cite{FS}).  In this note we will obtain some
very usable perturbation results for Weyl-Heisenberg frames with only
elementary assumptions by using the Wiener Amalgam space norm and by
adding continuity assumptions to the function $g$.  We will also give
examples to show that these results are best possible.

\section{Preliminary Results} \label{PR}

To simplify the notation a little, 
given a function $g\in \ltr$ and $a,b\in \mathbb R$, 
we define for $k\in \Bbb Z$ the function

$$
G_{k}(t) = \sum_{n\in \mathbb Z}g(t-na)\overline{g(t-na-k/b)}
$$
It is not difficult to prove that the series defining $G_k(t)$ converges
absolutely for a.e. $x$.
We will need the Weyl-Heisenberg Frame Identity
(see \cite{HW}, Theorem 4.1.5, or \cite{CL1} for a complete
treatment).
  
\begin{Thm} \label{WHFI} ({\bf WH-Frame Identity.})  If
$\sum_{n}|g(t-na)|^2 \le B$ a.e. and $f\in L^{2}(\Bbb R)$ is
bounded and compactly supported, then
\begin{eqnarray*}
&\bd  \sum_{n,m\in \Bbb Z}\ed &|<f,E_{mb}T_{na}g>|^{2} \\
& = &
b^{-1}\sum_{k\in \mathbb Z}
\int_{R}\overline{f(t)}f(t-k/b)\sum_{n}g(t-na)
\overline{g(t-na-k/b)}dt\\ 
& = & b^{-1}\int_{R}|f(t)|^{2}\sum_{n}|g(t-na)|^2 dt 
 +   b^{-1}\sum_{k\not=
0}\int_{R}\overline{f(t)}f(t-k/b)G_{k}(t)\ dt.
\end{eqnarray*}
\end{Thm}

We define the Wiener Amalgam space $W(L^\infty, {\ell}^1)$ as the
set of functions $g\in \ltr$ for which for some  $a>0$,
$$||g||_{W,a}:= \sum_n ||T_{na}g \cdot \chi_{[0,a[}||_\infty < \infty
.$$
It can be proved that if $||g||_{W,a}$ is finite for one value of $a$,
it
is automatically finite for all $a$. Furthermore, $||g||_{W,a}$ defines
a norm on $W(L^\infty, {\ell}^1)$.

\vspace{.2in}  We need some elementary facts about the Wiener
Amalgam space.  These can be found for example in \cite{HW}, Proposition
4.1.7
and the proof of Theorem 4.1.8.

\begin{Lem}\label{1.2}  Let $g\in W(L^\infty, \ell^1)$.

(1)  If $0<a\le b$ then $\|g\|_{W,b}\le 2\|g\|_{W,a}$.

(2)  $\|g\|_{W,a/2}\le 2\|g\|_{W,a}$.

(3)  Given functions $f,h\in W(L^{\infty},{\ell}^{1})$,
$$
\sum_{k}\left|\left|
\sum_{n}|T_{na}f||T_{na+k/b}h|\right|\right|_{\infty}\le
4 \|f\|_{W,a}\|h\|_{W,a}.
$$
\end{Lem}

The next result follows from the proof of Theorem 2.3 from
\cite{DSXW}.

\begin{Lem}\label{1.3}
For $g\in \ltr$ and bounded, compactly supported $f$, we have
\begin{eqnarray*}
\suz \int |\overline{f(t)} f(t - k/b)
\sun g(t-na) \overline{g(t -na-k/b)}|
dt \\ \le
\int|f(t)|^{2}\sum_{k\in \mathbb Z}|G_{k}(t)|\ dt.
\end{eqnarray*}
\end{Lem}

We also need the
 perturbation result of Christensen and Heil
\cite{CH} 

\begin{Thm} \label{1} If $(f_{i})$ is a frame  with frame bounds
$A,B$
 and there exists a constant $R\in [ 0, A[$
such that
for all $f\in H$,
$$
 \sum_{i}|<f,f_{i}-g_{i}>|^{2} \le
R\|f\|^{2},
$$
then $(g_i)$ is a frame with bounds
$A(1-\sqrt{\frac{R}{A}})^2, \ B(1+\sqrt{\frac{R}{B}})^2$.
\end{Thm}

\section{Perturbations}\label{Per}
\setcounter{equation}{0}

We start with a Proposition which contains the basic tool for
our first perturbation result. In light of theorem \ref{1} all we
really need to show is that the system $(h-g,a,b)$ has a finite upper
frame bound.  More specifically: 

\begin{Prop}\label{2.2}
Suppose $(g,a,b)$ is a WH-frame with frame bounds $A,B$ and let
$h\in \ltr$.  If there exists $R<A$ such that
\begin{eqnarray}
\sum_{k\in \Bbb Z}| \sum_{n\in \Bbb Z}(h-g)(t-na)
\overline{(h-g)(t-na -k/b)}|\le bR,
\ \label{3}
\text{a.e.},
\end{eqnarray}
then $(h,a,b)$ is a Weyl-Heisenberg  
frame for $H$ with frame bounds \\ $A(1-\sqrt{\frac{R}{A}})^2, 
B(1+\sqrt{\frac{R}{B}})^2$.  Moreover, if $(g,a,b)$ is a Riesz
basis for $L^{2}(\mathbb R )$, then $(h,a,b)$ is also a Riesz basis.
\end{Prop}

{\it Proof.} Let $f$ be bounded and compactly supported.
By the WH-frame Identity and Lemma  \ref{1.3}
we have:
\begin{eqnarray*}
&\bd \sum_{n,m\in \Bbb Z}\ed &|<f,E_{mb}T_{na}(h-g)>|^{2}\\
&=&\frac{1}{b}\sum_{k\in \Bbb Z}\int_{\Bbb R}\overline{f(t)}f(t-k/b)
\sum_{n\in \Bbb Z}(h-g)(t-na)\overline{(h-g)(t-na -k/b)} dt \\
&
\le &
\frac{1}{b} \int_{\Bbb R}|f(t)|^{2}
\sum_{k}|\sum_{n\in \Bbb Z}(h-g)(t-na)(h-g)\overline{(t-na -k/b)}|\ dt
\\
& \le & R\|f\|^{2}
\end{eqnarray*}
The set of bounded compactly supported functions is dense in $\ltr$, so
the above estimate actually holds for all functions $f\in \ltr$.
By Lemma \ref{1}, we have that
$(h,a,b)$ is a frame with the given bounds, and $(h,a,b)$ is a Riesz
basis if $(g,a,b)$ is a Riesz basis.
\qed \vspace{14pt}

In the paper of Jing \cite{J} there is a section concerning
perturbations of
Weyl-Heisenberg frames which at first glance appear to be similar
to our results.
For example,  in \cite{J} 
one of the main perturbation results for WH-frames is that 
if $(g,a,b)$ is a WH-frame and
$$ \left|\left| \sum_{k,n\in \Bbb Z} 
|(g-h)(\cdot -na- k/b)|^2 \right|\right|_\infty < bA,$$
then $(h,a,b)$ is also a frame. However, it should be observed that if
$ab$ is 
rational, this
condition is only satisfied if $g=h$ a.e., i.e., the result is not
useful 
in that case. 
Suppose namely that $(g-h)(x)\not= 0$.
Since there exist an infinite number of $n,k \in \mathbb Z$ such that
$na+\frac{k}{b}=0$, it follows that \\
$\sum_{k,n\in \mathbb Z} |(g-h)(x-na-\frac{k}{b}|^2=\infty$. However,
Proposition \ref{2.2} above applies for any value of $ab$.

We will now show that our perturbation result works whenever $g,h$ are
close
in the Wiener Amalgam norm.  Note that this result does not require $g$
to
be in the Wiener Amalgam space.

\begin{Thm}\label{2.4}
Suppose that $(g,a,b)$ is a WH-frame with frame bounds $A,B$. Let  
$h\in
\ltr$, and assume there exists $R<A$ such that
$$
\|g-h\|_{W,a}\le \sqrt{\frac{bR}{4}}.
$$
Then $(h,a,b)$ is a WH-frame with bounds $A(1-\sqrt{\frac{R}{A}})^2,
B(1+\sqrt{\frac{R}{B}})^2$.  Moreover, if $(g,a,b)$ is a Riesz basis 
for $L^{2}(\mathbb R )$, then
$(h,a,b)$ is also a Riesz basis.
\end{Thm}

{\it Proof.}
Using Lemma \ref{1.2}, we have 
\begin{eqnarray*}
& \bd \sum_{k\in \Bbb Z} \ed & | \sum_{n\in \Bbb Z}(h-g)(t-na)
\overline{(h-g)(t-na -k/b)}|  \\
&\le&  \sum_{k}\left|\left|  \sum_{n}|T_{na}(g-h)||T_{na+k/b}(g-h)| 
\right|\right|_\infty  \\
&\le&  4\|g-h\|_{W,a}\|g-h\|_{W,a} =
4\|g-h\|_{W,a}^{2} \le bR.  
\end{eqnarray*}

So the result follows from  Proposition \ref{2.2}.
\qed \vspace{14pt}

The condition $R<A$ in Proposition \ref{2.2} can not be relaxed.
To see this, fix ${\epsilon}>0$ and let
$$g = {\chi}_{[0,1]}+(1-{\epsilon}){\chi}_{[1,2]}, \ \
h = {\chi}_{[0,2]},
$$
\ $(g,1,1)$ is a Riesz basis for $\ltr$ with lower frame
bound
${\epsilon}^{2}$, since, for any finite sequence  of scalars
 $(a_{mn})_{m,n\in \mathbb Z}$ we have

\begin{eqnarray*}
\|\sum_{m,n\in \mathbb Z}a_{mn}E_{m}T_{n}g\| &\ge&
 \|\sum_{m,n\in \mathbb Z}a_{mn}E_{m}T_{n}{\chi}_{[0,1]}\| -  
\|\sum_{m,n\in \mathbb
Z}a_{mn}E_{m}T_{n}(1-{\epsilon}){\chi}_{[1,2]}\| \\  
&=& {\epsilon}\|\sum_{m,n\in \mathbb
Z}a_{mn}E_{m}T_{n}{\chi}_{[0,1]}\| = 
{\epsilon}\left ( \sum_{m,n\in \mathbb Z}|a_{mn}|^{2} \right )
^{1/2}. 
\end{eqnarray*}
 
Also,

$$
\sum_{k}|\sum_{n\in \mathbb Z}(h-g)(x-n)\overline{(h-g)(x-k-n)}| = 
{\epsilon}^{2}, \ \ \text{for
all}\ \ x.
$$

But $(h,1,1)$ is not a WH-frame.  The easiest way to check this is
to use the well known fact that $(h,1,1)$ is a WH-frame if and only
if it is a Riesz basis.  But,
$$
\|\sum_{k=0}^{2n-1}(-1)^{k}T_{k}h\| =
\|{\chi}_{[0,1]}-{\chi}_{[2n-1,2n]}\|
= \sqrt{2}.
$$
So $(T_{k}h)_{k\in \mathbb Z}$ is not a Riesz basic sequence in $L^{2}(
\mathbb R )$.

Let $(g,a,b)$ be a WH-frame.
It is an open question which Weyl-Heisenberg frames are equivalent to
compactly supported Weyl-Heisenberg frames.  Also, it is another
delicate
question when we can restrict the function $g$ to a compact subset of
$\mathbb R$ and still have a WH-frame for $L^{2}(\mathbb R )$.  
This question goes directly to the heart of applications where
compactly supported WH-frames are used.  Our next 
result shows that this is possible whenever $g\in
W(L^{\infty},{\ell}^{1})$.

\begin{Cor}
If $g\in W(L^{\infty},{\ell}^{1})$ and $(g,a,b)$ is a WH-frame,
then there is a natural number $N$ so that $({\chi}_{[-n,n]}g,a,b)$
is a WH-frame whenever $n \ge N$.
\end{Cor}

{\it Proof.}
We assume that $(g,a,b)$ is a WH-frame with frame bounds 
$A,B$.  Since $g\in W(L^{\infty},\ell^{1})$ we have
$$
\sum_{n\in \mathbb Z}\|{\chi}_{[n,n+1)}g\|_{\infty} <\infty,
$$
and so
$$
\lim_{n\rightarrow \infty}\sum_{|m|\ge n}\|{\chi}_{[m,m+1)}
g\|_{\infty}=
\lim_{n\rightarrow \infty}\|{\chi}_{[-n,n]}g-g\|_{W,a}= 0.
$$
Hence, there is an $N$ so that for all $n\ge N$ we have
that $\|{\chi}_{[n,n+1)}g\|_{\infty}\le 1$ and
$$
\|{\chi}_{[-N,N]}g - g\|_{W,1} < \sqrt{\frac{bA}{4}}.
$$
Now for $n\ge N$ we have
$$
\|{\chi}_{[-n,n]}g-g\|_{L^{2}(\mathbb R )}^{2} \le
\int_{|t|\ge N}|f(t)|^{2} \ dt
= 
\sum_{|k|\ge N}\int_{0}^{1}|f(t-k)|^{2} \ dt
$$
$$
\le \sum_{|k|\ge N}\|{\chi}_{[k,k+1)}g\|_{\infty}^{2}\le
\sum_{|k|\ge N}\|{\chi}_{[k,k+1)}g\|_{\infty}
= \|{\chi}_{[-N,N]}g - g\| < \sqrt{\frac{bA}{4}}.
$$ 
The Corollary now follows from Theorem \ref{2.4}. 
\qed \vspace{14pt}

Now we have a considerable strengthening of Proposition \ref{2.2} for
the case $a=b=1$.  

\begin{Thm}\label{11}
Let $(g,1,1)$ be a WH-frame with frame bounds $ A,B$.  
Let $h\in L^{2}(\mathbb R )$ and $0< {\lambda} < 1$ satisfy
$$
\sum_{n\in \mathbb Z}|(g-h)(x+n)| \le {\lambda}
\sqrt{A}\ \ \text{a.e.}
$$
Then $(h,1,1)$ is a WH-frame for $L^{2}(\mathbb R )$ with
frame bounds 
$$
(1-{\lambda})^{2}A\ \ \mbox{and} \ \ 
(1+{\lambda})^{2}B.
$$
\end{Thm}

{\it Proof.}
If $Z$ is the Zak transform, we have
$$
|Z(g)(x,y)-Z(h)(x,y)|= |\sum_{n\in \mathbb Z}g(x+n)e^{2{\pi}iny}-
\sum_{n\in \mathbb Z}h(x-n)e^{2{\pi}iny}| 
$$
$$
\le \sum_{n\in \mathbb Z}|(g-h)(x+n)| \le {\lambda}\sqrt{A}
\le {\lambda}|Zg(x,y)|.
$$
It follows that,
$$
(1-{\lambda})\sqrt{A}\le (1-{\lambda})|Z(g)(x,y)|\le |Z(h)(x,y)|
$$
$$
\le
(1+{\lambda})|Z(g)(x,y)|\le (1+{\lambda})\sqrt{B}.
$$
So $(h,1,1)$ is a Weyl-Heisenberg frame for 
$L^{2}(\mathbb R)$ with the stated frame bounds (see \cite{HW}, Theorem
4.3.3).  \qed \vspace{14pt}

It is easily seen that we can not allow $\lambda=1$ in the inequality in 
Theorem \ref{11}.  For example, if $g={\chi}_{[0,1]}, \
h={\chi}_{[0,2]}$
then $(g,1,1)$ is an orthonormal basis for $L^{2}(\mathbb R )$ and
as we saw earlier, $(h,1,1)$ is not a frame.  But,
$$
\sum_{n\in \mathbb Z}|(g-h)(x+n)| = 1\ \ \text{a.e.}
$$

We might hope for an even sharper result with the inequality in
Theorem \ref{11} changed to
$$
\sum_{n\in \mathbb Z}|(g-h)(x+n)|^{2}\le {\lambda}A^{\alpha},
$$
for some $0<{\alpha}\le 1$.
Unfortunately, this fails.  For example, let
$$
g = {\chi}_{[0,1]}
$$
and
$$
h = \frac{1}{2}{\chi}_{[0,2]}.
$$
Then $(h,1,1)$ is not a frame (since 
$(T_{n}h)_{n\in \mathbb Z}$ is not a Riesz basic sequence) while
$(g,1,1)$ is an orthonormal basis for $L^{2}(\mathbb R )$ (and
so $A=B=1$).
Finally, 
$$
\sum_{n\in \mathbb Z}|(g-h)(t+n)|^{2} = \left ( \frac{1}{2} \right )
^{2}
+ \left ( \frac{1}{2} \right ) ^{2} = \frac{1}{2}=\frac{1}{2}A^{\alpha}.
$$

Now let $(g,a,b)$ be a frame
and we will look at perturbations of the modulation and
translation parameters to see when we can still be guaranteed to
have a WH-frame.  The main problem here is that we may not be
able to change
$a$ or $b$ by any arbitrarily small amount and still get a frame.  This
follows from a result of Feichtinger and Janssen \cite{FJ}. They show
that 
there is a  function
$g\in L^{2}(\Bbb R)$ so that $(g,a,b)$ has a finite upper frame
bound only
when a and b are rational.  Therefore, no matter how close
$(a^{'},b^{'})$ is
to $(a,b)$, we still may not have a frame.  The next
technical difficulty occurs if $a=b=1$.  If $(g,1,1)$ is a WH-frame,
then
we can never have a general result of the form: $|a'-a|<{\epsilon}$
implies
$(g,a',1)$ is a frame since if $a'>1$ then  $(g,a',1)$
is never complete.  Despite these strong limitations, we can obtain some
satisfactory perturbation results which will guarantee that if the
translation parameters are close enough then we will have a frame for
all
small $b$.  In this result, as well as the rest of the results in this
section, the price we pay for being able to perturb in one
parameter is that the other parameter may change drastically.

\begin{Thm}\label{2.7}
Let $g\in W(L^{\infty}, \ell^1)$ with $(g,a,b)$ a WH-frame 
with frame bounds $A,B$ and
let $0<R< bA $. There is an 
$0< {\epsilon}\le \frac{a}{2}$ and $b_{0} = b_{0}({\epsilon})$ so
that whenever 
$|a-a'|<{\epsilon}$ and
$$
\sum_{n}|g(t-na)-g(t-na')|^{2}\le R ,\ \ \text{a.e.},
$$
then $(g,a',b')$ is a WH-frame
whenever $0<b' <b_{0}$.
\end{Thm}

{\it Proof.}
If $(g,a,b)$ generates a WH-frame 
with frame bounds $A,B$ then (see Heil and Walnut \cite{HW},
the proof of Proposition 4.1.4, page 649) 
$$
bA\le  \sum_{n\in \mathbb Z}|g(t-na)|^{2} \le bB,\ \ \text{a.e.}
$$
Using the (reverse) triangle inequality we have

\begin{eqnarray*}
\sqrt{bA}-\sqrt{R} &\le &
\left ( \sum_{n}|g(t-na)|^{2}\right ) ^{1/2} -
\left ( \sum_{n}|g(t-na)-g(t-na')|^{2}\right ) ^{1/2} \\
&\le & \left ( \sum_{n}|g(t-na')|^{2}\right ) ^{1/2} \\
& \le &\left ( \sum_{n}|g(t-na)|^{2}\right ) ^{1/2} +
\left ( \sum_{n}|g(t-na)-g(t-na')|^{2}\right ) ^{1/2} \\
& \le &  \sqrt{bA}+\sqrt{R}  \text{ a.e. }
\end{eqnarray*}

For the rest, we borrow an argument from \cite{HW} (the proof
of Theorem 4.1.8).
 Fix $0< {\epsilon}\le a/2$ satisfying
$$
{\delta}=: 32{\epsilon}\|g\|_{W,a}+16{\epsilon}^{2}<
[ \sqrt{bA}-\sqrt{R}]^{2}.  
$$
Now let $N$ be so large that
$$
\sum_{|n|\ge N}\|g\cdot {\chi}_{[an,a(n+1))}\|_{\infty} < \epsilon .
$$
Let $g_{0}= g\cdot {\chi}_{[-aN,aN]}$ and $g_{1} = g-g_{0}$, so 
that $\|g_{1}\|_{W,a}< {\epsilon}$.  Now if 
$$
b \le \frac{1}{4aN} = b_{0}
$$
then (with $G_k^\prime(t):=\sum_{n}T_{na'}g(x)\cdot 
T_{na'+k/b'}\overline{g(x)}$)

\begin{eqnarray*}
&\bd \sum_{k\not= 0} \ed & \|G_k^\prime (t)\|_{\infty} =
\sum_{k\not= 0}\left \| \sum_{n}T_{na'}g\cdot T_{na'+k/b'}\overline{g}
\right \| _{\infty} \\
 &\le & \sum_{k\not= 0}\left \| \sum_{n}|T_{na'}g||T_{na'+k/b'}g \right
\|  _{\infty} \\
 &=& \sum_{k\not= 0}\left \| \sum_{n}|T_{na'}g_{0} + 
T_{na'}g_{1}||T_{na'+k/b'}g_{0}+T_{na'+k/b'}g_{1}|\right \|_{\infty}\\ 
 &\le&  \sum_{k\not= 0}\left \|
\sum_{n}|T_{na'}g_{0}||T_{na'+k/b'}g_{0}|  \right \| _{\infty} + 
\sum_{k\not= 0}\left \| \sum_{n}|T_{na'}g_{0}||T_{na'+k/b'}g_{1}| 
\right \| _{\infty}+\\
&+& \sum_{k\not= 0}\left \| \sum_{n}|T_{na'}g_{1}||T_{na'+k/b'}g_{0}| 
\right \| _{\infty} + \sum_{k\not= 0}\left \|
\sum_{n}|T_{na'}g_{1}||T_{na'+k/b'}g_{1}|
\right \| _{\infty} \\
&\le&  0 + 8\|g_{0}\|_{W,a'}\|g_{1}\|_{W,a'} + 4\|g_{1}\|^{2}_{W,a'}.
\end{eqnarray*}
Now, since $\frac{a}{2}\le a' \le 2a$, we can continue our inequality
using Lemma \ref{1.2}, (1) and (2) to get:

\begin{eqnarray*}
\sum_{k\not= 0}\|G_{k}^{'}(t)\|_{\infty} & = & \\
\sum_{k\not= 0}\left \| \sum_{n}T_{na^{'}}g\cdot
T_{na'+k/b'}\overline{g}
\right \| _{\infty} & \le & 8\|g_{0}\|_{W,a'}\|g_{1}\|_{W,a'} + 
4\|g_{1}\|^{2}_{W,a'} \\
&\le & 32\|g_{0}\|_{W,a}\|g_{1}\|_{W,a} + 16\|g_{1}\|^{2}_{W,a}\\
&\le & 32{\epsilon}\|g\|_{W,a}+16{\epsilon}^{2} = \delta.
\end{eqnarray*}
It follows by Lemma \ref{1.3} that if $|a-a'|<\epsilon $ and $0<b' \le
b_{0}$ 
then
for all bounded, compactly supported functions $f\in L^{2}(\mathbb R )$
we have
\begin{eqnarray*}
\frac{1}{b^\prime}& \bd | \sum_{k\not= 0} \ed &\int_{\mathbb
R}\overline{f(t)}f(t-k/b')G_k^\prime(t)\ dt| \\
&\le& \frac{1}{b^\prime}\sum_{k\not= 0}\|G_{k}^\prime(t)\|_{\infty}
\int_{\mathbb R} 
|f(t)|^{2}\ dt \le \frac{1}{b^\prime}{\delta}\|f\|^{2}.
\end{eqnarray*} 
Also, from the first part of the proof, for all $f$ as above we have,

$$
\frac{1}{b^\prime}(\sqrt{bA}-\sqrt{R})^{2}\|f\|^{2}\le 
\frac{1}{b^\prime}\int_{\mathbb R}|f(t)|^{2}\sum_{n\in \mathbb 
Z}|g(t-na')|^{2} \
dt
\le \frac{1}{b^\prime}(\sqrt{bA}+\sqrt{R})^{2}\|f\|^{2}.
$$
Finally, the WH-Frame Identity yields,

\begin {eqnarray*} \sum_{m,n\in \mathbb Z}  |\langle
f,E_{mb'}T_{na'}g\rangle|^{2}  
&=& \frac{1}{b^\prime}\int_{\mathbb R}|f(t)|^{2}\sum_{n\in 
\mathbb Z}|g(t-na')|^{2}\ dt \\
&+& \frac{1}{b^\prime}\sum_{k\not= 0}\int_{\mathbb
R}\overline{f(t)} f(t-k/b')G_n^\prime(t)\ dt.
\end{eqnarray*}
Putting this altogether we have that
$$
\frac{1}{b^\prime}[(\sqrt{bA}-\sqrt{R})^{2}-{\delta}]\|f\|^{2} \le
\sum_{m,n\in \mathbb Z}|\langle f,E_{mb'}T_{na'}g\rangle |^{2} \le
\frac{1}{b^\prime}[(\sqrt{bA}+\sqrt{R})^{2}+{\delta}]\|f\|^{2}.
$$
Since this inequality holds for all bounded compactly supported
functions
$f\in L^{2}(\mathbb R )$, it holds for all $f\in L^{2}(\mathbb R )$,
which
completes the proof. 
\qed \vspace{14pt}

A general setting where the conditions of Theorem \ref{2.7} will
hold is when $g$ is continuous.  This is just enough to offset
the Feichtinger-Janssen example \cite{FJ}.

\begin{Cor}
If $g\in W(L^{\infty},\ell^1)$ is continuous and $(g,a,b)$ is a frame,
then
there is a ${\delta}>0$ and a $b_{0}>0$ so that
 $(g,a',b')$ is a WH-frame  whenever
$$
|a-a'|< \delta ,
$$
 and $0<b'<b_{0}$.
\end{Cor}

{\it Proof.}
We just need to verify that the conditions of Theorem \ref{2.7}
hold.  Fix $R<bA$.  Since $f\in W(L^{\infty},\ell^{1})$, we can
choose a natural number $n_{0}$ so that
$$
\|(1-{\chi}_{[a(-n_{0}+1),an_{0}]})g\|_{W,a}< \frac{R}{3}.
$$
Since $g$ is continuous on the compact set
$[-an_{0},a(n_{0}+1)]$, it
is uniformly continuous there.  In particular, there is a ${\delta}>0$
so that if $x,y\in [-an_{0},a(n_{0}+1)]$ then
$$
|x-y|\le {\delta},\ \ \Rightarrow \ \ |g(x)-g(y)|^{2}<  
\frac{R}{3(2n_{0}+2)}.
$$
Let ${\epsilon}=\frac{\delta}{n_{0}}$.
Now, if $|a-a'|<\epsilon$ and  then for all $-n_{0}\le n \le
n_{0}-1$ we have $$
|(t-na)-(t-na')|= |n||a-a'|<|n|{\epsilon}= \frac{|n|}{n_{0}}{\delta}\le
{\delta}. $$
Hence, for $t\in [0,a]$,
$$
|g(t-na)-g(t-na')|^{2}< \frac{R}{3(2n_{0}+2)},
$$
It follows that
\begin {eqnarray*}
& \bd \sum_{n}\ed & |g(t-na)-g(t-na')|^{2} \\
&=&\sum_{n=-n_{0}}^{n_{0}+1}|g(t-na)-g(t-na')|^{2}  
+\sum_{{n<-n_{0}}\atop {n>n_{0}+1}}|g(t-na)-g(t-na')|^{2}\\
&\le&
(2n_{0}+2)\frac{R}{3(2n_{0}+2)} +
 2\|(1-{\chi}_{[a(-n_{0}+1),an_{0}]})g\|_{W,a}\\
&<& \frac{R}{3}+ \frac{2R}{3} = R.
\end{eqnarray*}
The Corollary now follows by Theorem \ref{2.7}.
\qed \vspace{14pt}

We now have immediately the corresponding result for compactly supported
functions.

\begin{Cor}
If $(g,a,b)$ is a WH-frame where $g$ is compactly supported and
continuous, then there is a ${\delta}>0$ and a $b_{0}>0$ so
that $(g,a',b')$ is a WH-frame whenever
$$
|a-a'|< \delta ,
$$
 and $0<b'<b_{0}$.
\end{Cor}

Continuity is necessary in the preceding results.  A trivial example
occurs if we consider $({\chi}_{[0,1]},1,1)$ since no matter how close
we have $a'$ to $a$, if $a< a'$, we cannot have a frame for any $b$
since
a necessary condition for $(g,a,b)$ to form a WH-frame 
with frame bounds $A,B$ is that $
bA\le  \sum_{n\in \mathbb Z}|g(t-na)|^{2} \le bB,\ \ \text{a.e.}
$
In light of this, it is more natural to ask for
$(g,a',b')$ to form a frame for $0<a-a' < {\epsilon}$, and all
small $b'$.  But again, the above results will fail without the
assumption of continuity.  For example, we can let
$$
E_{1} = [0,1-\frac{1}{16}),
$$
$$
\ \ \ E_{2} = \cup_{n=2}^{\infty}
[1-\frac{1}{2^{2n}},1-\frac{1}{2^{2n+1}}),
$$
and
$$ 
E_{3} = \cup_{n=2}^{\infty}
[2-\frac{1}{2^{2n+1}},2-\frac{1}{2^{2(n+1)}}).
$$
Let $F=E_{1}\cup E_{2}\cup E_{3}$ and $g = {\chi}_{F}$.  Then it
is immediate that $(g,1,1)$ is an orthonormal basis for $L^{2}(\mathbb R
)$.
Now, if 
$$
1-\frac{1}{2^{2n+1}}< a' < 1-\frac{1}{2^{2(n+1)}}
$$
then for all 
$$
1-\frac{1}{2^{2n+1}}< t \le a',
$$
we have that $g(t) = 0$, and for $n\ge 1$, $t-na'<0$ so 
$g(t-na') = 0$.  Also, for $n\ge 2$ we have that
$2< t+na'$ and so $g(t+na') = 0$.  Finally, for $n=1$ we have that
$$
 2-\frac{1}{2^{2n}} = 1-\frac{1}{2^{2n+1}} + 1-\frac{1}{2^{2n+1}}
\le t+a' < 1-\frac{1}{2^{2(n+1)}}+ 1-\frac{1}{2^{2(n+1)}}=
2-\frac{1}{2^{2n+1}}.
$$
Hence, $g(t+a') = 0$.  It follows that
$$
\sum_{n\in \mathbb Z}|g(t-na')|^{2} = 0, \ \ \text{for all}\ \ 
1-\frac{1}{2^{2n+1}}< t \le a'.
$$
In particular, $(g,a',b)$ is not a frame for all $0< b$.  It follows
that,
given any ${\epsilon}>0$, there is an interval of points $a'$ with
$0< a-a' < \epsilon$ so that $(g,a',b)$ is not a frame for all $0< b$.

\end{document}